\documentclass[letterpaper, 10 pt, conference]{ieeeconf}  

\IEEEoverridecommandlockouts                              
\usepackage{graphics} 
\usepackage{epsfig} 
\usepackage{amsmath} 
\usepackage{amssymb}  
\usepackage{xcolor}
\usepackage{psfrag}
\usepackage{overpic}

\usepackage{amsmath, amssymb}

\newcommand{\bq}{\mathbf{q}}
\newcommand{\bp}{\mathbf{p}}

\newcommand{\bu}{\mathbf{u}}

\newcommand{\bx}{\mathbf{x}}

\newcommand{\bB}{\mathbf{B}}
\newcommand{\bC}{\mathbf{C}}

\newcommand{\bQ}{\mathbf{Q}}
\newcommand{\bR}{\mathbf{R}}

\newcommand{\bU}{\mathbf{U}}
\newcommand{\bV}{\mathbf{V}}

\newcommand{\bX}{\mathbf{X}}

\newcommand{\bxi}{\boldsymbol{\xi}}
\newcommand{\bXi}{\boldsymbol{\Xi}}

\newcommand{\bSigma}{\boldsymbol{\Sigma}}
\newcommand{\bTheta}{\boldsymbol{\Theta}}

\newcommand{\bGamma}{\boldsymbol{\Gamma}}

\def\stackrel#1#2{\mathrel{\mathop{#2}\limits^{#1}}}

\title{\LARGE \bf
Discovering conservation laws from data for control
}

\author{Eurika Kaiser$^{1}$ , J. Nathan Kutz$^{2}$  and Steven L. Brunton$^{3}$ 
\thanks{$^{1}$ Eurika Kaiser is with the Department of Mechanical Engineering, University of Washington, Stevens Way, United States
        {\tt\small eurika@uw.edu}}%
\thanks{$^{2}$ Nathan Kutz is with the Department of Applied Mathematics, University of Washington, Stevens Way, United States
	{\tt\small kutz@uw.edu}}
\thanks{$^{3}$ Steven Brunton is with the Department of Mechanical Engineering, University of Washington, Stevens Way, United States
	{\tt\small sbrunton@uw.edu}}%
}

\begin{document}

\maketitle
\thispagestyle{empty}
\pagestyle{empty}

\begin{abstract}
Conserved quantities, i.e. constants of motion, are critical for characterizing many dynamical systems in science and engineering.  These quantities are related to underlying symmetries and they provide fundamental knowledge about physical laws, describe the evolution of the system, and enable system reduction.  In this work, we formulate a data-driven architecture for discovering conserved quantities based on Koopman theory. The Koopman operator has emerged as a principled linear embedding of nonlinear dynamics, and its eigenfunctions establish intrinsic coordinates along which the dynamics behave linearly. Interestingly, eigenfunctions of the Koopman operator associated with vanishing eigenvalues correspond to conserved quantities of the underlying system. 
In this paper, we show that these invariants may be identified with data-driven regression and power series expansions, based on the infinitesimal generator of the Koopman operator.  
We further establish a connection between the Koopman framework, conserved quantities, and the Lie-Poisson bracket.  
This data-driven method for discovering conserved quantities is demonstrated on the three-dimensional rigid body equations, where we simultaneously discover the total energy and angular momentum and use these intrinsic coordinates to develop a model predictive controller to track a given reference value.  
\end{abstract}

\section{INTRODUCTION}
Conservation laws lie at the heart of the physical sciences.  
The seminal work of Emmy Noether has had a profound impact on our understanding of physics, establishing a deep connection between conservation laws and the underlying symmetries of the system~\cite{noether1918invariante}. 
For example, many Hamiltonian systems are completely integrable, meaning that they are entirely described by the integrals of motion that remain constant along the flow of the dynamics~\cite{Abraham1978book,Marsden:MS,arnol2013mathematical}.  
These integrals of motion establish an intrinsic description that does not depend on the specific coordinate system, which is useful for symmetry reduction~\cite{Marsden:MS}, simulation~\cite{Yoshida1990pla,Marsden2001dmvi}, and control~\cite{Kaiser2017arxiv}.  
However, conservation laws have been historically elusive, and are difficult to discover, even with access to the governing equations.  
Recently, connections have been established between conserved quantities and modern Koopman operator theory~\cite{Koopman1931pnas,Mezic2005nd}, which naturally lends itself to data-driven characterizations~\cite{Mezic2017book,Kaiser2017arxiv}.  
In this work, we leverage modern data-driven techniques, including sparse optimization~\cite{Brunton2016pnas,Kaiser2017arxiv}, to discover and control conserved quantities in the Koopman operator framework.  

Advances in data-driven modeling and control are opening up new possibilities for the automated discovery of governing equations and conservation laws~\cite{Langley1981bacon,Schmidt2009science,schulte2010discovery,Brunton2016pnas,Rudy2017sciadv,Kaiser2017arxiv}. 
The BACON discovery algorithm~\cite{Langley1981bacon}, developed in 1981, uses expert input to design experiments and has learned many basic laws including Snell's law of refraction, Black's specific heat law, and Joule's law.  
In the seminal work of Schmidt and Lipson~\cite{Schmidt2009science}, genetic programming is used to automatically discover the conserved energy in a double pendulum from experimental data.   
More recently, sparsity promoting optimization has been used to discover equations of motion~\cite{Brunton2016pnas} and conserved quantities, such as the Hamiltonian, in a range of systems~\cite{Kaiser2017arxiv}.  
Sparse optimization provides a particularly parsimonious description, identifying the fewest terms required to explain the data.  
The resulting models are both interpretable and prevent overfitting.  
Conserved quantities are also naturally expressed in the Koopman operator perspective~\cite{Koopman1931pnas,Mezic2005nd,Mezic2017book,Mezic2017arxiv}, which provides an intrinsic description of nonlinear dynamical systems in terms of an infinite-dimensional linear operator acting on the space of measurements.  
In particular, conserved quantities are eigenfunctions of the Koopman operator with zero eigenvalue.  
The Koopman operator is amenable to data-driven characterization~\cite{Schmid2010jfm,Rowley2009jfm,Williams2015jnls,Kutz2016book} and provides an embedding where it is possible to design controllers for nonlinear systems in a linear control framework~\cite{Kaiser2017arxiv,Korda2016arxiv,Peitz2017arxiv}.  
For conservative systems, the conserved quantities themselves provide a particularly useful coordinate system in which to enact control~\cite{Kaiser2017arxiv}. 

In this work, we leverage data-driven optimization to discover and control conserved quantities in dynamical systems following the subsequent steps:
\begin{enumerate}
	\item Data-driven identification of multiple conserved quantities in the null-space of the infinitesimal generator equation. We assume here full access to the state.
	\item Identification of the control term in the equation for the time evolution of the identified conserved quantities, which constitutes the model for the controller.
	\item Design of, e.g., a model predictive controller (MPC)~\cite{garcia1989model,camacho2013model}, in these intrinsic coordinates.
\end{enumerate}
We also establish a connection between Koopman theory and the Lie-Poisson bracket, which may be used to characterize the space of conserved quantities, provides a consistency check for candidate invariants, and can be used to estimate the vector field.  
This work extends previous work~\cite{Kaiser2017arxiv} by (i) discovering multiple conserved quantities, (ii) estimating the control term from data, and (iii) introducing usable connections to the Lie-Poisson bracket formalism.
We demonstrate these concepts on the three-dimensional rigid body equations.
However, the resulting theory and methods are general, and may broadly be applied to develop parsimonious models and effective controllers for space mission design~\cite{koon2000heteroclinic,koon2008dynamical}, inviscid flows~\cite{Marsden:MS},  and various mechanical systems such as robots~\cite{Bloch1996arma}, to name a few. 
\section{BACKGROUND}
\subsection{Koopman spectral theory}
\label{Sec:KoopmanSpectralTheory}
We consider a nonlinear, deterministic dynamical system of the form
\begin{equation}\label{Eq:Dynamics}
\frac{\mathrm d}{\mathrm dt}\bx(t) =  {\boldsymbol{\it f}}(\bx)
\end{equation}
with the state of the system $\bx\in\mathcal{M}$, where $\mathcal{M}$ is a differentiable manifold, often given by $\mathcal{M}=\mathbb{R}^n$.  
In discrete time, the dynamics are given by 
%
$\bx(t_0+t) = \mathbf{F}_t(\bx(t_0)),$ 
%
where $\mathbf{F}_t$ may be the flow map $\mathbf{F}_t(\bx(t_0)) = \bx(t_0)+\int_{t_0}^{t_0+t} {\boldsymbol{\it f}}(\bx(\tau))\,d\tau$ induced by the dynamics in \eqref{Eq:Dynamics}.
%
%
In 1931, B. O. Koopman~\cite{Koopman1931pnas} introduced the operator theoretic perspective for dynamics, showing that there exists an infinite-dimensional linear operator $\mathcal{K}_t$ that acts to advance the Hilbert space of all scalar measurement functions $g(\bx)$ along the flow $\mathbf{F}_t$ of the dynamics:  
\begin{equation}
\mathcal{K}_t g(\bx(t_0)) = g(\mathbf{F}_t (\bx(t_0))) = g(\bx(t_0+t)).\label{Eq:KoopmanDiscrete2}
\end{equation}
Thus, the Koopman operator is a composition operator, $\mathcal{K}_tg = g\circ\mathbf{F}_t$.
For smooth dynamics, there is a continuous system
\begin{equation}\label{Eqn:KoopmanGenerator}
\frac{\mathrm d}{\mathrm dt} g = \mathcal{A}g,
\end{equation}
where $\mathcal{A}:={\bf f}\cdot\nabla$ is the infinitesimal generator of the family of Koopman operators $\mathcal{K}_t$, parameterized by time $t$.  

The linearity of the Koopman operator is appealing. However, its infinite-dimensional nature poses issues for representation and computation.  
Instead of capturing the evolution of all measurement functions in a Hilbert space, applied Koopman analysis approximates the evolution on a subspace spanned by a finite set of measurement functions, e.g. through eigen-observables or intrinsic measurement functions~\cite{Williams2015jnls,Brunton2016plosone}.  
The spectral decomposition of the Koopman operator thus facilitates reduced system representation and provides crucial insights into the dynamics, e.g. invariants, stability, eigenmeasures, etc.  

A Koopman eigenfunction $\varphi(\bx)$ corresponding to eigenvalue $\lambda$ satisfies
%
$\lambda\varphi(\bx) = \varphi(\mathbf{F}_t(\bx))$
%
and in continuous-time satisfies
\begin{equation}
\frac{\mathrm d}{\mathrm dt}\varphi(\bx) = \lambda \varphi(\bx).\label{Eq:KoopmanEfun}
\end{equation}
Obtaining Koopman eigenfunctions analytically or from data is a central challenge in modern dynamical systems.  
Importantly, discovering these eigenfunctions enables globally linear representations of strongly nonlinear systems.  

In recent work~\cite{Kaiser2017arxiv}, we showed that eigenfunctions can be identified with data-driven regression and power series expansions, based on the infinitesimal generator of the Koopman operator.
%
%
Combined with \eqref{Eq:KoopmanEfun}, an eigenfunction $\varphi(\bx)$ for a specific eigenvalue $\lambda$ satisfies the linear partial differential equation (PDE):
\begin{equation}
{\nabla\varphi(\bx)\cdot\mathbf{f}(\bx)= \lambda\varphi(\bx)}.\label{Eq:KoopmanPDE}
\end{equation}
This formulation assumes that the eigenfunctions are smooth.  
A data-driven procedure to identify eigenfunctions from \eqref{Eq:KoopmanPDE} can be found in~\cite{Kaiser2017arxiv} and is briefly outlined in Sec.~\ref{Sec:IdentifyFromData}.

\subsection{Conservation laws and symmetries}
Koopman operator theory does not assume a specific structure of the system~\eqref{Eq:Dynamics}. 
However, in the present work we focus on Hamiltonian systems, and more broadly on measure-preserving systems. 

For an autonomous, finite-dimensional Hamiltonian systems, the propagator has a symplectic structure given by ${\boldsymbol{\it f}} = [ \partial H/\partial\bp\;\;-\partial H/\partial\bq]^T$, where $H$ is the Hamiltonian. The canonical form of the dynamics is then given by
\begin{equation}\label{Eqn:EOM_Hamiltonian}
	\frac{\mathrm d}{\mathrm dt}\bq = \frac{\partial H}{\partial\bp},\quad
	\frac{\mathrm d}{\mathrm dt}\bp = -\frac{\partial H}{\partial\bq}, 
\end{equation}
where $\bq$ is the generalized coordinate and $\bp$ is the conjugate momentum.
Hamiltonian systems are a prominent class of measure-preserving systems, for which phase space volume is preserved.  

Conservation laws provide critical information about the dynamics and the integrability of the system, and are related to underlying symmetries.
Specifically, a system with $n$ degrees of freedom is fully integrable if there exist $n$ functionally independent first integrals, a subset of constants of motion.
While many problems do not admit full integrability, even knowing some constants of motion provides invaluable information about the physical behavior of these systems. 

Through the celebrated Noether theorem~\cite{noether1918invariante}, intrinsic symmetries are related to conserved quantities in the Lagrangian framework, with Lagrangian $L(\bq,\bp)$.  
We assume no explicit time dependency for simplicity.
A coordinate $q^i$ is considered {\it cyclic} if $\partial L/\partial q^i = 0$. If the total time derivative $\mathrm{d} p_i/\mathrm{dt} = 0$ vanishes, where $p_i = \partial L/\partial \dot{q}^i$ is the conjugate momentum to $q^i$, then this coordinate is conserved. 
These ideas are generalized in Noether's theorem.
Thus, if the coordinate $q^i$ is cyclic, this variable does not affect the Lagrangian and the Lagrangian is invariant with respect to translations of that variable, so that the Lagrangian possesses a symmetry~\cite{Marsden:MS}. 
Moreover, conservation of the conjugate momentum follows, hence all coordinates are not linearly independent and a constraint equation exists.  
This constraint can be used to eliminate the cyclic coordinate from the equations, resulting in a symmetry reduction~\cite{Marsden:MS}.  
An example is the total energy, which is related to symmetry in time, i.e.\ translations in time do not have any effect, and result in the conservation of energy.

A function $C_i:\mathbb{R}^{2n}\rightarrow\mathbb{R}^n$ is a constant of motion if it Poisson-commutes with the Hamiltonian function $H$:
\begin{equation}\label{Eqn:PoissonBracket_CQ}
\{ C_i,H\} = \sum\limits_{j=1}^n \left( \frac{\partial C_i}{\partial q_j} \frac{\partial H}{\partial p_j} -  \frac{\partial C_i}{\partial p_j} \frac{\partial H}{\partial q_j}\right) = 0.
\end{equation}
It is also said that $C_i$ and $H$ are in mutual involution.
For full integrability, the $n$ constants of motion must all be in mutual involution, i.e.\ $\{C_i,C_j\}=0\,\forall i,j$.
A function $C$ that Poisson-commutes with every function $F$, i.e.\  $\{ C,F\} = 0$ for any $F$,  
is called Casimir function or invariant.
Casimir functions are crucial for identifying constraints of the system, integrability, system reduction, and establishing stability criteria via the Energy-Casimir method~\cite{Marsden:MS}.


\section{CONNECTIONS BETWEEN CONSERVED QUANTITIES AND THE KOOPMAN OPERATOR}
\label{Sec:CQ_and_KO}
Here we consider a scalar observable $g(\bx)$ of the state, parametrized in terms of the positions $\bq$ and momenta $\bp$.  
The total derivative with respect to time is:
\begin{equation}
\frac{\mathrm d}{\mathrm dt}g = \frac{\partial g}{\partial \bx}\frac{\partial \bx}{\partial t} 
= \frac{\partial g}{\partial \bq} \frac{\partial \bq}{\partial t}  
+ \frac{\partial g}{\partial \bp} \frac{\partial \bp}{\partial t} 
+ \frac{\partial g}{\partial t}.
\end{equation}
The last term vanishes, if we assume that there is no explicit time dependency.
We can identify the remaining terms with the (canonical) Poisson bracket: 
\begin{equation}
\frac{\mathrm d}{\mathrm dt}g 
= \frac{\partial g}{\partial \bq} \frac{\partial H}{\partial \bp}  
- \frac{\partial g}{\partial \bp} \frac{\partial H}{\partial \bq}
= \{ g, H \}.
\end{equation}
Thus, the evolution of any measurement function with no explicit time dependency in a Hamiltonian system is given by its Poisson bracket with the Hamiltonian.
In addition, the dynamics of the measurement function is governed by the Koopman operator~\eqref{Eqn:KoopmanGenerator}, 
$\frac{\mathrm d}{\mathrm dt}g 
= \{ g, H \}
= \mathcal{A} g$.
For smooth eigenfunctions of the generator, we then obtain
\begin{equation}
\frac{\mathrm d}{\mathrm dt}\varphi
= \{ \varphi, H \}
= \lambda \varphi.
\end{equation}
Thus, the Hamiltonian is an eigenfunction of the Koopman operator associated with $\lambda=0$, which follows from
$\frac{\mathrm d}{\mathrm dt}H
= \{ H, H \}=\lambda H = 0$
for $\varphi = H$ and
where the Poisson bracket vanishes due to the symplectic structure.  
Further, any other conserved quantity $C_i$ is also a Koopman eigenfunction with $\lambda=0$:
\begin{equation}\label{TimeDerivative_ConservedQuantity}
\frac{\mathrm d}{\mathrm dt}C_i
= \{ C_i, H \}
= \lambda C_i  = 0, 
\end{equation}
which follows from the Poisson-commutativity \eqref{Eqn:PoissonBracket_CQ}.
Thus, Koopman eigenfunctions associated with the eigenspace at $\lambda=0$ remain constant along the flow of the dynamics, and are therefore conserved quantities.   

There may not always exist smooth constants of motion; however, non-smooth, discontinuous ones are more general and may exist and we refer to~\cite{mezic1999chaos,Mezic2004physicad} for further reading and~\cite{mezic2003cdc} for their use in the context of control. 
 



\section{IDENTIFYING CONSERVATION LAWS FROM DATA}
\label{Sec:IdentifyFromData}
We propose two approaches to identify conserved quantities from data using sparse regression:
(A) based on the Koopman reduced order nonlinear identification for control (KRONIC) framework for identifying Koopman eigenfunctions~\cite{Kaiser2017arxiv}, and
(B) using the Lie-Poisson bracket formalism. 
In both cases, the conserved quantities lie in the null space of a constructed data matrix.
More generally, both approaches can be used to learn other Koopman eigenfunctions.

\subsection{Based on the PDE for the infinitesimal generator of the Koopman operator}\label{Sec:KRONIC}
Building on the sparse identification of nonlinear dynamics (SINDy) framework~\cite{Brunton2016pnas}, Koopman eigenfunctions for a particular value of $\lambda$ can be identified using the PDE~\eqref{Eq:KoopmanPDE}.  

First, a dictionary of candidate functions is chosen:
\begin{equation}
\bTheta(\bx) = \begin{bmatrix} \theta_1(\bx) &\theta_2(\bx) & \cdots & \theta_p(\bx)\end{bmatrix}.
\end{equation}
The dictionary $\bTheta$ must be large enough and carefully chosen so that a Koopman eigenfunction may be well approximated:
\begin{equation}
\varphi(\bx) \approx \sum_{k=1}^p\theta_k(\bx)\xi_k = \bTheta(\bx)\bxi.
\end{equation}
%
Given data $\bX = [\bx_1\, \bx_2\, \cdots\, \bx_m]$, where $\bx_i:=\bx(t_i)$ denotes the sampled measurement of the system~\eqref{Eq:Dynamics},  
the time derivative $\dot{\bX}=[\dot{\bx}_1\, \dot{\bx}_2\, \cdots\, \dot{\bx}_m]$ can be approximated numerically from $\bx(t)$ if not measured directly~\cite{Brunton2016pnas}. 
The total variation derivative~\cite{Chartrand2011isrnam} is recommended for noise-corrupted measurements.
The dictionary evaluated on the data becomes:
%
\begin{equation}\label{Eq:Theta}
\bTheta(\bX) = \begin{bmatrix} \theta_1(\bX^T) & \theta_2(\bX^T) & \cdots & \theta_p(\bX^T)\end{bmatrix}.
\end{equation}
Moreover, we can define a library of directional derivatives, representing the possible terms in $\nabla\varphi(\bx)\cdot\mathbf{f}(\bx)$ from~\eqref{Eq:KoopmanPDE}:  $\bGamma(\bx,\dot{\bx})=[\nabla\theta_1(\bx)\cdot\dot{\bx}\; \nabla\theta_2(\bx)\cdot\dot{\bx}\; \cdots\; \nabla\theta_p(\bx)\cdot\dot{\bx} ]$.  It is then possible to construct $\bGamma$ from data, 
$\bGamma(\bX,\dot{\bX}) =
\scriptsize
\begin{bmatrix} \nabla\theta_1(\bX^T)\cdot\dot{\bX} & \nabla\theta_2(\bX^T)\cdot\dot{\bX} & \cdots & \nabla\theta_p(\bX^T)\cdot\dot{\bX}\end{bmatrix}
$.

%

For a specific eigenvalue $\lambda$, the Koopman PDE in \eqref{Eq:KoopmanPDE} evaluated on data yields:
\begin{equation}\label{Eq:SparseKoopman}
\left(\lambda\bTheta(\bX) - \bGamma(\bX,\dot{\bX})\right)\bxi = \mathbf{0}.
\end{equation}
The formulation in \eqref{Eq:SparseKoopman} is implicit, so that $\bxi$ will be in the null space of ${\lambda\bTheta(\bX)-\bGamma(\bX,\dot{\bX})}$.  
The null space can be determined via singular value decomposition (SVD). Specifically,
the right null space is spanned by the right singular vectors of $\lambda\bTheta(\bX)-\bGamma(\bX,\dot{\bX}) = \bU \bSigma\bV^{*}$ (i.e., columns of $\bV$) corresponding to zero-valued singular values.  
The few active terms in an eigenfunction can be identified by imposing a sparsity constraint on the coefficient vector, e.g. using alternating direction methods~\cite{Qu2014}.

Constants of motion belong to eigenfunctions of the eigenspace at $\lambda=0$ (see Sec.~\ref{Sec:CQ_and_KO}). Thus it is sufficient to search for a solution vector in the null space:
\begin{equation}
 - \bGamma(\bX,\dot{\bX})\bxi = \mathbf{0}.
\end{equation}
The set of vectors that are mapped to zero by $\bGamma(\bX,\dot{\bX})$ is denoted by the kernel, $\mathrm{ker}(\bGamma)$. 
The dimension of the kernel, $\mathrm{dim}(\mathrm{ker}(\bGamma))$, corresponds to the number of integrals of motion.
Thus, any linear combination of the vectors in the kernel constitutes a constant of motion.
Determining the eigenspace at zero does not mean there exist a finite-dimensional representation of the Koopman operator; in contrast, a system, e.g. a frictionless pendulum, may have conserved quantities such as the Hamiltonian, but exhibits otherwise only a continuous spectrum of the Koopman operator.


\subsection{Based on the Lie-Poisson bracket formalism}
The Lie-Poisson bracket formalism can be employed to identify conserved quantities using ideas developed in~\cite{Kaiser2017arxiv}, with canonical Hamiltonian systems being a special case. 
The following is derived for a Lie-Poisson bracket associated with non-canonical Hamiltonian systems, which arise in fields such as dynamical systems~\cite{guckenheimer_holmes} and fluid dynamics~\cite{Thiffeault2000physicad}.
It is straightforward to formulate this for the Poisson bracket of systems with a canonical structure.
In particular, this formulation can be used (1) to identify other conserved quantities based on a known conserved quantity, and (b) as a consistency check for the identified subspace of conserved quantities.

The dynamics of a non-canonical Hamiltonian system can be written in terms of the Lie-Poisson bracket:
\begin{equation}\label{Eqn:LiePoissonBracket_Dynamics}
\dot{x}_i
= \{ x_i, H \} = \langle \bx,\nabla x_i\times\nabla H\rangle,\;i=1,\ldots,n,
\end{equation}
where the Lie-Poisson bracket is defined as:
\begin{equation}\label{Eqn:LiePoissonBracket}
\{ F,G \} (\boldsymbol{\bx}) := \langle-\boldsymbol{\bx}, \nabla F\times\nabla G\rangle = -\bx\cdot(\nabla F\times\nabla G)
\end{equation}
for two scalar functions $F$ and $G$.
For later reference, we define the partial derivative of a scalar function $F$ in terms of its candidate basis:
\begin{eqnarray}\label{Eqn:BasisDerivative}
\partial_i F({\bf x}) 
:=\boldsymbol{\Theta}_i({\bf x}) \boldsymbol{\xi} =  \begin{bmatrix}
\partial_i\theta_1(\bx),\partial_i\theta_2(\bx),\ldots,\partial_i\theta_p(\bx)
\end{bmatrix} \boldsymbol{\xi},
\end{eqnarray}
where $\partial_i:= \frac{\partial}{\partial x_i}\forall i$ and $\nabla \boldsymbol{\Theta}({\bf x}) \boldsymbol{\xi} = [\boldsymbol{\Theta}_1({\bf x}) \boldsymbol{\xi}, \ldots,\boldsymbol{\Theta}_n({\bf x}) \boldsymbol{\xi}]$ follows.

%

We assume that the Hamiltonian function has been previously identified using KRONIC, i.e.\ $H = \boldsymbol{\Theta}(\boldsymbol{\bx})\boldsymbol{\xi}$, then the vector field can be estimated:
\begin{equation}\label{Eqn:LiePoissonBracket_VectorField}
\dot{x}_i
= f_i({\bf x}) = \{ x_i, \boldsymbol{\Theta}\boldsymbol{\xi} \} = \langle \bx,\nabla x_i\times\nabla \boldsymbol{\Theta}\boldsymbol{\xi}\rangle,\;\forall i.
\end{equation}

Further, we can use the knowledge of the Hamiltonian to identify another constant of motion $C$, which shall be represented as:
\begin{equation}\label{Eqn:H_ansatz}
C = \boldsymbol{\Upsilon}(\boldsymbol{\bx})\boldsymbol{\eta}, 
\end{equation}
where 
$\boldsymbol{\Upsilon}$ is another library:
\begin{equation}
\boldsymbol{\Upsilon}(\bx) = 
\begin{bmatrix}
\upsilon_1(\bx),\upsilon_2(\bx),\ldots,\upsilon_r(\bx)
\end{bmatrix}.
\end{equation}
The number and type of candidate functions can generally be different, , i.e. $r{\not = }p$ and $\theta_i{\not = }v_j\,\forall i,j$.\\[0.5ex]

\noindent\textbf{Derivation in three degrees of freedom.} In order to identify $C$, we will insert the ansatzes~\eqref{Eqn:H_ansatz} 
 into~\eqref{TimeDerivative_ConservedQuantity} and make use of the derivative~\eqref{Eqn:BasisDerivative}: 
\begin{eqnarray}
\dot{C} (\bx) 
&=& -\bx\cdot 
\begin{pmatrix}
\partial_3 H\partial_2 C - \partial_2 H\partial_3 C\\
\partial_1 H\partial_3 C - \partial_3 H\partial_1 C\\
\partial_2 H\partial_1 C - \partial_1 H\partial_2 C\\
\end{pmatrix}\notag\\
&=& x_1(\partial_2 H\partial_3 C - \partial_3 H\partial_2 C) \notag\\
&& \hspace{-0.25cm}+ x_2(\partial_3 H\partial_1 C - \partial_1 H\partial_3 C)\notag\\
&& \hspace{-0.25cm}+ x_3(\partial_1 H\partial_2 C - \partial_2 H\partial_1 C) \stackrel{!}{=} 0,
\end{eqnarray}
and with the libraries:
\begin{eqnarray}
0
&=& x_1(
\boldsymbol{\Theta}_2\boldsymbol{\xi} \boldsymbol{\Upsilon}_3 \boldsymbol{\eta} 
- \boldsymbol{\Theta}_3\boldsymbol{\xi} \boldsymbol{\Upsilon}_2 \boldsymbol{\eta})\notag\\
&& \hspace{-0.25cm}+ x_2(
\boldsymbol{\Theta}_3\boldsymbol{\xi} \boldsymbol{\Upsilon}_1 \boldsymbol{\eta} 
- \boldsymbol{\Theta}_1\boldsymbol{\xi} \boldsymbol{\Upsilon}_3 \boldsymbol{\eta})\notag\\
&& \hspace{-0.25cm}+ x_3(
\boldsymbol{\Theta}_1\boldsymbol{\xi} \boldsymbol{\Upsilon}_2 \boldsymbol{\eta} 
- \boldsymbol{\Theta}_2\boldsymbol{\xi} \boldsymbol{\Upsilon}_1 \boldsymbol{\eta})\notag\\
&=& \left[x_1(
\boldsymbol{\Theta}_2\boldsymbol{\xi} \boldsymbol{\Upsilon}_3
- \boldsymbol{\Theta}_3\boldsymbol{\xi} \boldsymbol{\Upsilon}_2)\right. \notag\\
&& \hspace{-0.25cm}+ x_2(
\boldsymbol{\Theta}_3\boldsymbol{\xi} \boldsymbol{\Upsilon}_1
- \boldsymbol{\Theta}_1\boldsymbol{\xi} \boldsymbol{\Upsilon}_3)\notag\\
&&\underbrace{ \left.\hspace{-0.25cm}+ x_3(
	\boldsymbol{\Theta}_1\boldsymbol{\xi} \boldsymbol{\Upsilon}_2
	- \boldsymbol{\Theta}_2\boldsymbol{\xi} \boldsymbol{\Upsilon}_1)\right]}_{\boldsymbol{D}} \boldsymbol{\eta}. 
\end{eqnarray}
%
Then, we are searching for a sparse vector $\boldsymbol{\eta}$ in the null space of $\boldsymbol{D}$ to identify the active components for $C$ in~\eqref{Eqn:H_ansatz}.
In addition, evaluation of the Lie-Poisson bracket is useful as a consistency check for identified constants of motion using KRONIC, as constants of motion must be in mutual involution. 



\section{CONTROL IN INTRINSIC COORDINATES}\label{Sec:ControlDesign}
The dynamics of an affine-in-control system are given by
\begin{equation}\label{Eq:DynamicsControl}
\frac{\mathrm d}{\mathrm dt}\bx(t) =  {\boldsymbol{\it f}}(\bx) + \bB\bu
\end{equation}
with continuously differentiable dynamics $\boldsymbol{\it f}:\mathbb{R}^n\rightarrow\mathbb{R}^n$, control matrix $\bB\in\mathbb{R}^{n\times q}$, and multi-channel input $\bu\in\mathbb{R}^q$.

Applying the chain rule to a Koopman eigenfunction yields
\begin{subequations}\label{Eq:ControlEfun}
	\begin{align}
	\frac{\mathrm d}{\mathrm dt}\varphi(\bx) 
	&=  \nabla\varphi(\bx)\cdot {\boldsymbol{\it f}}(\bx) + \nabla\varphi(\bx)\cdot\bB\bu\\
	&= \lambda\varphi(\bx) + \nabla\varphi(\bx)\cdot\bB\bu.
\end{align}
\end{subequations}
This formulation of the control problem has also been considered in~\cite{Surana2016cdc} in the context of observer design and observability.
For conserved quantities $C_i$, this becomes
\begin{subequations}\label{Eq:ControlCQ}
	\begin{align}
	\frac{\mathrm d}{\mathrm dt}C_i(\bx) 
	= \nabla C_i(\bx)\cdot\bB\bu = \bB_{c}(\bx) \bu,
	\end{align}
\end{subequations}
where the control matrix $\bB_{c}(\bx):=\nabla C_i(\bx)\cdot\bB$ is generally nonlinear.
For known eigenfunctions, such as conserved quantities, 
control may be designed and applied in these {\it intrinsic} coordinates.
In terms of the KRONIC formalism, the dynamics in~\eqref{Eq:ControlCQ} become
\begin{subequations}\label{Eqn:RBM:CQ_CONTROL_EQNS}
	\begin{align}
		\frac{d}{dt}{\bC}(\bx)   &= \boldsymbol{\Theta}_{\bx}(\bx)\boldsymbol{\Xi}\cdot \bB \bu,
	\end{align}
\end{subequations}
where $\bTheta_{\bx}(\bx):=\begin{bmatrix}
\nabla\theta_1(\bx) & \nabla\theta_2(\bx) & \ldots & \nabla\theta_p(\bx)
\end{bmatrix}$ denotes the gradient of the dictionary of candidate functions.
Then, model-based control strategies can be employed, e.g.\ by solving a state-dependent Riccati equation to find a sub-optimal controller~\cite{cimen2008ifac} or by using model predictive control~\cite{garcia1989model,camacho2013model}.
We refer to~\cite{Kaiser2017arxiv} for a non-affine formulation.\\[0.5ex]

\noindent\textbf{Discovering the control matrix ${\bf B}$ from data.} In general, the control matrix ${\bf B}$ may be unknown and it is of interest to discover it from data.
Building on the identification of Koopman eigenfunctions in Sec.~\ref{Sec:KRONIC}, we can then estimate ${\bf B}$ from~\eqref{Eq:ControlEfun} or~\eqref{Eqn:RBM:CQ_CONTROL_EQNS} from forced data consisting of sampled pairs $\{{\bf x}_i,{\bf u}_i \}_{i=1}^M$ of the state and the actuation input.
It can be shown that \eqref{Eq:ControlEfun} can be re-arranged as
\begin{eqnarray}\label{Eqn:IdentifyB}
\dot{\varphi}
=  \lambda \varphi + \nabla \varphi\cdot 
\begin{bmatrix}
- & {\bf b}_1 & -\\
& \vdots & \\
- & {\bf b}_n & -
\end{bmatrix}{\bf u}
=  \lambda \varphi +\nabla \varphi\otimes  {\bf u}^T 
\begin{bmatrix}
{\bf b}_1^T\\
\vdots\\
{\bf b}_n^T
\end{bmatrix}
\end{eqnarray}
where ${\bf b}_i$ represents the $i$th row of ${\bf B}$ and
${\bf A}\otimes {\bf b}^T$ is defined as $ (\begin{smallmatrix}
a_{11} & a_{12}\\
a_{21} & a_{22}
\end{smallmatrix})\otimes (\begin{smallmatrix}
b_1 & b_2
\end{smallmatrix}) := (\begin{smallmatrix}
a_{11}[b_1\,b_2] & a_{12}[b_1\,b_2]\\
a_{21}[b_1\,b_2] & a_{22}[b_1\,b_2]
\end{smallmatrix}) = (\begin{smallmatrix}
a_{11}b_1 & a_{11}b_2 & a_{12}b_1 &a_{12}b_2\\
a_{21}b_1 & a_{21}b_2& a_{22}b_1 & a_{22}b_2
\end{smallmatrix})$.
The matrix ${\bf B}\in\mathbb{R}^{n\times q}$ is now represented as a vector in $\mathbb{R}^{nq\times 1}$.
Assuming the eigenfunction $\varphi({\bf x})=\boldsymbol{\Theta}({\bf x}){\boldsymbol{\xi}_i}$ and associated eigenvalue $\lambda$ have been discovered through the KRONIC architecture using unforced data, 
then~\eqref{Eqn:IdentifyB} can be written in terms of the libraries as
\begin{eqnarray}\label{Eqn:IdentifyBdata} 
\left[\boldsymbol{\Gamma}({\bf x}, \dot{\bf x}) - \lambda\boldsymbol{\Theta}({\bf x})\right]\boldsymbol{\xi}
= \nabla \boldsymbol{\Theta}({\bf x})\boldsymbol{\xi}\otimes  {\bf u}^T 
\begin{bmatrix}
{\bf b}_1^T\\
\vdots\\
{\bf b}_n^T
\end{bmatrix}.
\end{eqnarray}
While $\boldsymbol{\xi}$ has been discovered using KRONIC on unforced data, 
$\boldsymbol{\Theta}({\bf x})$ and $\boldsymbol{\Gamma}({\bf x}, \dot{\bf x})$ are here evaluated on data collected from the forced system; 
then $\dot{\varphi} = \lambda \varphi + \nabla \varphi({\bf x})\cdot ({\bf f}({\bf x})+{\bf B}{\bf u}) = \boldsymbol{\Gamma}({\bf x}, \dot{\bf x})\boldsymbol{\xi}$.
This can now be used to solve for ${\bf b}_i$, $i=1,\ldots,n$ by computing the Moore-Penrose pseudoinverse of  $\nabla \boldsymbol{\Theta}({\bf x})\boldsymbol{\xi}\otimes  {\bf u}^T$.
For the special case that an eigenfunction corresponds to a conserved quantity, $\lambda$ is set to zero in \eqref{Eqn:IdentifyBdata}.


\section{EXAMPLE: RIGID BODY SYSTEM}


We consider the rotation of a rigid body in a coordinate frame attached to its center of mass and rotating with the body.
The states of the system are the angular momentum $\Pi_i = I_i\omega_i$, where $I_i$ are the principal moments of inertia for $i=1,2,3$.  In these coordinates, as opposed to angular velocity $\omega_i$ or Euler angles, the Hamiltonian structure of the Euler equations become apparent in terms of the Lie-Poisson bracket.

%
%
The change in angular momentum subject to an externally applied torque $\boldsymbol{\tau} = (\tau_1,\tau_2,\tau_3)$ are given by the forced Euler equations
%
\begin{equation}\label{Eqn:EulerEquationsWithControl}
	\dot{\boldsymbol{\Pi}}  = 
	\begin{bmatrix}
	\frac{I_2-I_3}{I_3 I_2}\Pi_2\Pi_3\\[1ex]
	\frac{I_3-I_1}{I_1 I_3}\Pi_3\Pi_1\\[1ex]
	\frac{I_1-I_2}{I_2 I_1}\Pi_1\Pi_2
	\end{bmatrix}
	+ \bB
	\begin{bmatrix}
	\tau_1\\
	\tau_2\\
	\tau_3
	\end{bmatrix}
	= {\it {\bf f}}(\boldsymbol{\Pi}) + \bB\boldsymbol{\tau},
\end{equation}
where $\bB$ is a three-dimensional identity matrix and $[I_1, I_2, I_3]= [1,1/2,1/3]$.
We aim to identify conserved quantities for the case of free rotation without external torque, i.e. $\boldsymbol{\tau} = [0,0,0]^T$. 

\begin{figure}[tb]
	\centering
	\begin{overpic}[height = 3.95cm,trim = {20 22 16.5 10}, clip=true]{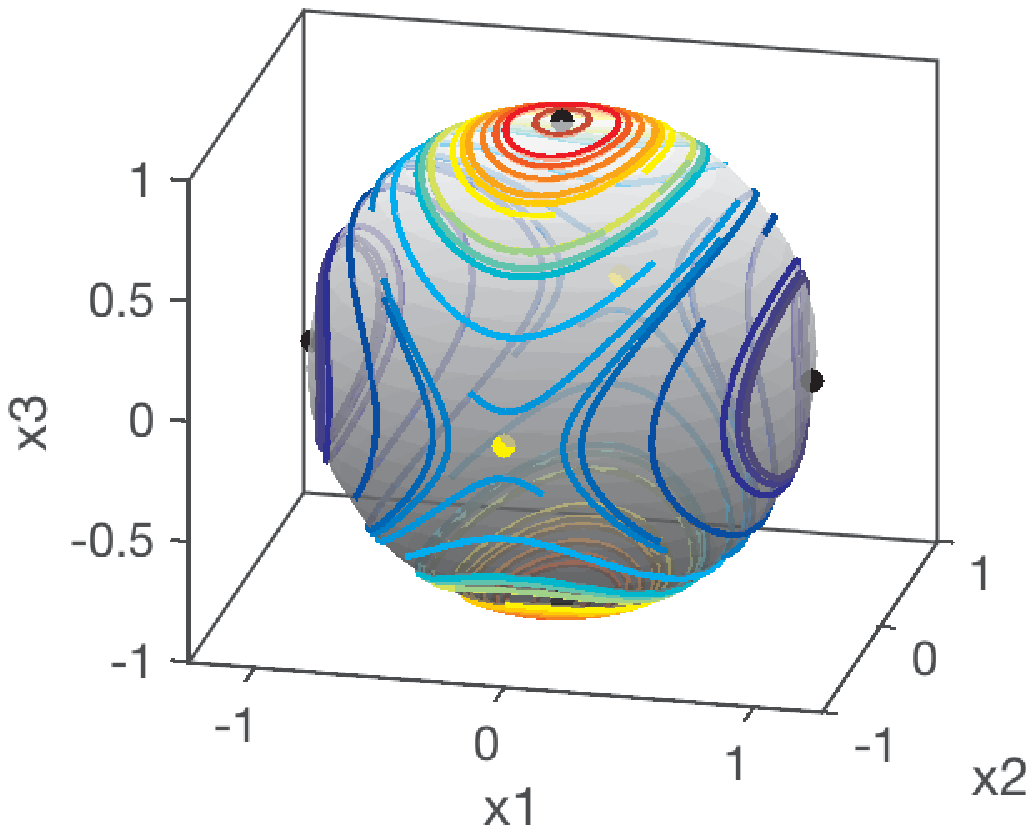}
		\put(2,-2.5){(a)}
		\put(-5,38){$\Pi_3$}
		\put(44,-4.5){$\Pi_1$}
		\put(97,20){$\Pi_2$}
	\end{overpic}
	\hfill
	\begin{overpic}[height = 3.95cm,trim = {20 22 16.5 10}, clip=true]{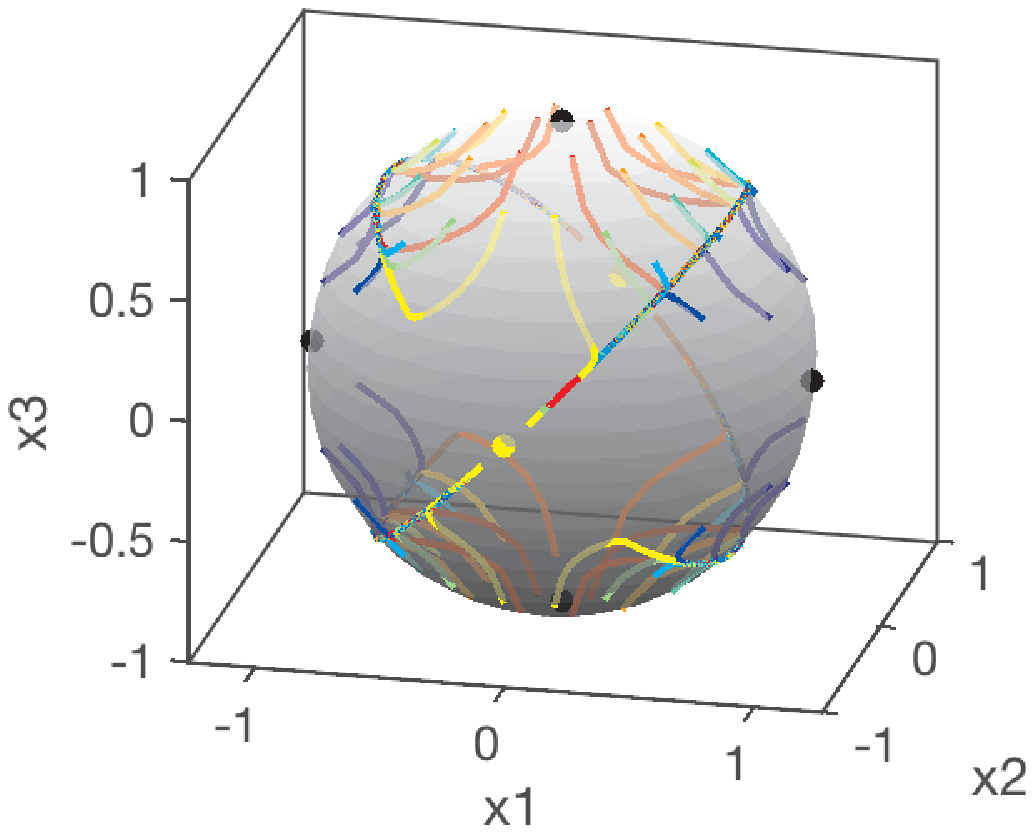}
		\put(2,-2.5){(b)}
	\end{overpic}\\[2ex]
	\psfrag{Time}{Time}
	\psfrag{H,L}{\hspace{-1mm}$H,\,L$}
	\psfrag{0}{$0$}
	\psfrag{0.5}{\hspace{-1mm}$0.5$}
	\psfrag{1}{$1$}
	\psfrag{1.5}{\hspace{-1mm}$1.5$}
	\psfrag{2}{$2$}
	\hspace{-0.cm}\begin{overpic}[height = 3.1cm, trim = {0 0 15 0}, clip=true]{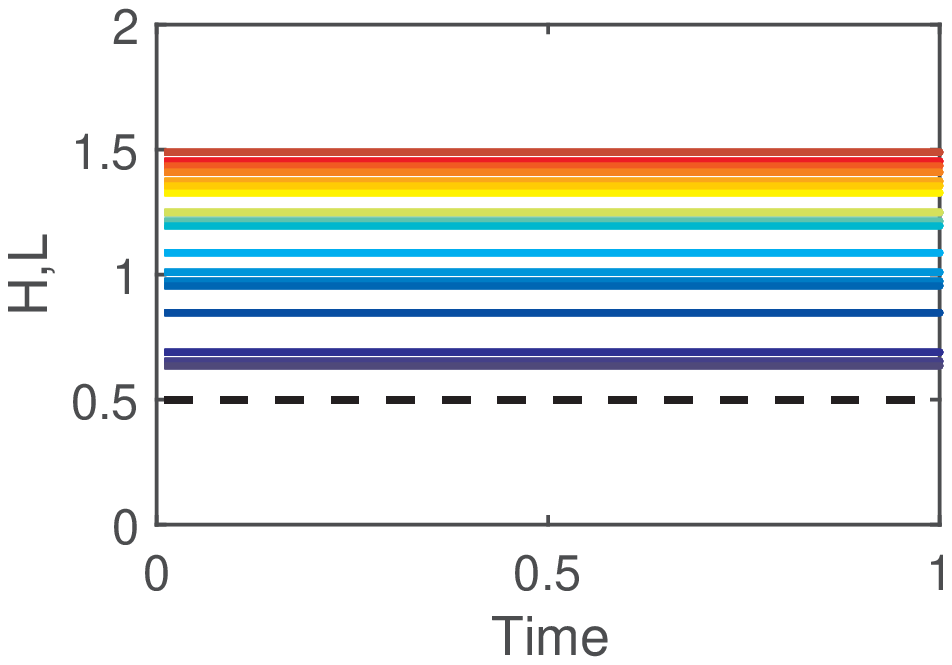}
		\put(2,0){(c)}
	\end{overpic}	
	\hfill
	\psfrag{Time}{Time}
	\psfrag{H,L}{}
	\psfrag{0}{}
	\psfrag{0.5}{}
	\psfrag{1}{}
	\psfrag{1.5}{}
	\psfrag{2}{}
	\begin{overpic}[height = 3.1cm, trim = {43 0 15 0}, clip=true]{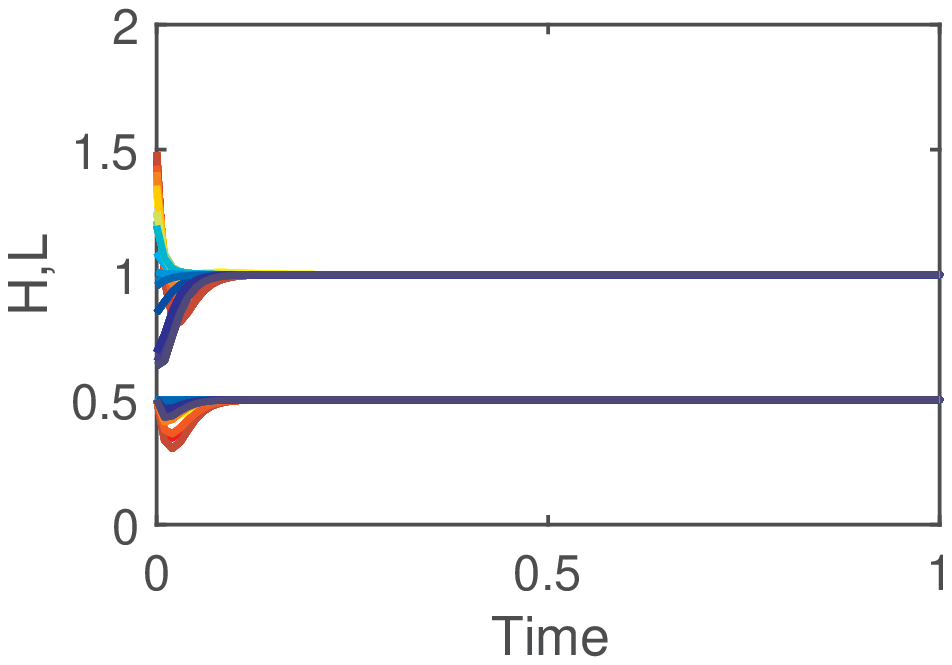}
		\put(-16,0){(d)}	
	\end{overpic}	
	\caption{Rigid body motion: 
		(a) unforced trajectories on the angular momentum sphere $L=const.$,
		(b) with KRONIC using the identified subspace of conserved quantities ${\bf C}$ and $\hat{\bf B}$ from data, 
		and (c-d) the evolution of the kinetic energy $H$ (colored lines) and angular momentum $L$ ($--$ in (c) and colored lines in (d))
		corresponding to each case and initial condition.}
	\label{Fig:RBM:ControlResults}
	\vspace{-3ex}
\end{figure}

%
%
Two quadratic first integrals of the system are the angular momentum  and the total kinetic energy (Hamiltonian):
\begin{equation}
\footnotesize
L(\boldsymbol\Pi) = \frac{1}{2}\left(\Pi_1^2 + \Pi_2^2 + \Pi_3^2\right),\; H(\boldsymbol\Pi) = \frac{1}{2}\left( \frac{\Pi_1^2}{I_1} + \frac{\Pi_2^2}{I_2} + \frac{\Pi_2^2}{I_2} \right).
\end{equation}
These confine the motion to the intersection of energy ellipsoids with momentum spheres, shown in Fig.~\ref{Fig:RBM:ControlResults} (a). 
Although the Euler equations do not exhibit a symplectic structure, 
there exists a simple Hamiltonian form in terms of the Lie-Poisson bracket~\eqref{Eqn:LiePoissonBracket}.
For instance, the conservation of angular momentum can be shown:
\begin{align}\label{Eqn:LiePoissonBracket_AngMom}
\dot{L} (\boldsymbol{\Pi}) 
&= \{ L,H \} (\boldsymbol{\Pi})
= -\boldsymbol{\Pi} \cdot (\nabla L \times \nabla H) \notag\\
&= -\boldsymbol{\Pi}\cdot (\boldsymbol{\Pi}\times \nabla H ) = 0,
\end{align}
and the Euler equations~\eqref{Eqn:EulerEquationsWithControl} can be recovered through
\begin{equation}
\dot{\Pi}_i = \{ \Pi_i,H \} (\boldsymbol{\Pi}) = -\boldsymbol{\Pi}\cdot(\nabla \Pi_i\times\nabla H),\; i=1,2,3.
\end{equation}

\subsection{Identification of conserved quantities and control matrix}
\begin{figure}[tb]
	\centering
	\psfrag{sigma}{$\sigma_k$}
	\psfrag{index}{$k$}
	\psfrag{0}{$0$}
	\psfrag{10}{$10$}
	\psfrag{20}{$20$}
	\psfrag{-5}{\hspace{-1mm}$-5$}
	\psfrag{5}{$5$}
	\includegraphics[width=0.4\textwidth]{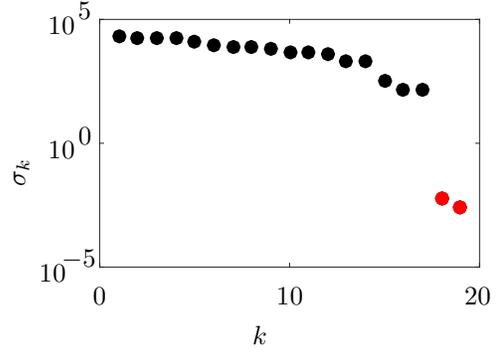}
	\caption{Identification of the null space ($\color{red}\bullet$) through SVD.}
	\label{Fig:RGM:SingularValues}
\end{figure}

Conserved quantities are identified from data for the unforced system~\eqref{Eqn:EulerEquationsWithControl} with $\boldsymbol{\tau} = [0,0,0]^T$.
An ensemble of short-time trajectories is collected over time $t\in\{0,10\}$ with time step $\Delta t = 0.01$.  Using the KRONIC formalism~\eqref{Eq:SparseKoopman}, conserved quantities are found in the nullspace of $\lambda\bTheta(\bX)-\bGamma(\bX,\dot{\bX})$ for $\lambda=0$. 
A polynomial basis up to order $p=3$ is employed to construct the library of candidate functions.
The null space is defined by the last two right singular vectors in the matrix $\bV$ of the SVD corresponding to nearly zero singular values (see Fig.~\ref{Fig:RGM:SingularValues}). 
Thus, the null space is given by $\hat{\boldsymbol{\Xi}} =[\boldsymbol{\xi}_1,\boldsymbol{\xi}_2]= [\begin{smallmatrix}0.58 & 0.58 & 0.58\\ 0.71 & 0 & 0.71  \end{smallmatrix}]^T$ depicting only the non-zero entries (all others are  $\approx \mathcal{O}(10^{-16})$).
%
The angular momentum is given by the first column of $\hat{\bXi}$, while the total energy is given by a linear combination of the columns since columns of $\bXi$ span an orthogonal subspace.
Thus, the nullspace yields conserved quantities up to a scaling. 
Since any linear combination of conserved quantities and any conserved quantity shifted by a constant is conserved, the identified null space represents the (infinite-dimensional) space of the conserved quantities.
The Lie-Poisson bracket formulation can be used as a consistency check that the identified null space really spans the space of conserved quantities. It can be shown that $\{ C_1,C_2\} = \{ \boldsymbol{\Theta}\boldsymbol{\xi}_1,\boldsymbol{\Theta}\boldsymbol{\xi}_2\} = \mathcal{O}(10^{-4})$.

The identified conserved quantities can then be used to discover the control matrix ${\bf B}$ using data sampled from the forced system.
In particular, an ensemble of trajectories is collected with simultaneous forcing in the three control inputs, $\boldsymbol{\tau}= [(0.5+\sin(40t))^3, 0.5+\sin(10t), \sin(20t)]^T$.
The matrix ${\bf B}$ is determined by solving for its entries in~\eqref{Eqn:IdentifyBdata}, where $\boldsymbol{\Theta}({\bf x})$ and $\boldsymbol{\Gamma}({\bf x},\dot{\bf x})$  are now evaluated on the data of the forced system and using the identified $\hat{\boldsymbol{\Xi}}$. The entries of the estimated ${\bf B}$ have an error of $\mathcal{O}(10^{-3})$.


\subsection{Control in intrinsic coordinates}
Having identified this subspace of conserved quantities $\boldsymbol{C}$ and an estimate for ${\bf B}$, it is now possible to design a controller for the system in these coordinates . Considering the Euler equations with externally applied torque as input~\eqref{Eqn:EulerEquationsWithControl},
%
%
the dynamics for the $C_i$ are given by~\eqref{Eq:ControlCQ}, which is used as predictive model on which basis the control input is determined.  
Specifically, we employ model predictive control to minimize the cost function:
\begin{equation}\label{Eqn:CostFunction}
	J = \int_{0}^{T}\,(\bC(\boldsymbol{\Pi})-\bC^{*})^T\bQ(\bC(\boldsymbol{\Pi})-\bC^{*}) + \boldsymbol{\tau}^T\bR\boldsymbol{\tau}\,dt
\end{equation}
with state and input weight matrices, $\bQ\in\mathbb{R}^{n\times n}$  and $\bR\in\mathbb{R}^{q\times q}$, respectively. Both matrices are symmetric and fulfill $\bQ > 0$ and $\bR\geq 0$.
In particular, we choose 
$\bQ = (\begin{smallmatrix}
2 & 0\\
0 & 2
\end{smallmatrix})$ and 
$\bR = \mathrm{diag}(0.001,0.001,0.001)$. 
Typically, a reference state is given, here $\boldsymbol{\Pi}^{*} = (0,1,0)$, on which the conserved quantities are then evaluated, i.e. $C_i^{*}:=C_i(\boldsymbol{\Pi}^{*})$, which constitutes the reference in~\eqref{Eqn:CostFunction}.
The control horizon is $T=10\Delta t$ and the dynamics~\eqref{Eqn:EulerEquationsWithControl} are integrated with time step $\Delta t/10$. The control is updated every $10$ integration steps and kept constant at intermediate steps.

The resulting controller is demonstrated on an ensemble of $N=114$ initial conditions and compared with the unforced dynamics (see Fig.~\ref{Fig:RBM:ControlResults}).
Initial conditions are chosen on the angular momentum sphere $L=1$.
Reference states are displayed as yellow circles and correspond to the unstable states to be stabilized. 
Note that the controller cannot distinguish between $\boldsymbol{\Pi}^{*}$ and $-\boldsymbol{\Pi}^{*}$ as these have the same energy and angular momentum.
The controller successfully drives all trajectories to the reference state.


\section{CONCLUSIONS AND FUTURE WORKS}
In summary, we have demonstrated a data-driven sparse optimization architecture for the simultaneous discovery and control of multiple conserved quantities, which we apply to the rigid body equations.
These conserved quantities are critical to reduce the equations of motion, for integrability, to characterize phase space organization, and to define intrinsic coordinates.
In addition, we establish connections between the Koopman operator, conserved quantities, and the Lie-Poisson bracket. 
Conservation laws are eigenfunctions of the Koopman operator corresponding to the zero eigenvalue, and thus characterizing conserved quantities benefits from the tremendous progress in Koopman analysis using modern data-driven and machine learning techniques. 
Model verification represents a critical, but often neglected step in the model identification process, and the Lie-Poisson bracket provides a consistency check for candidate discovered conserved quantities.
The resulting conserved quantities establish a natural coordinate system, in which control may be developed, facilitating the data-driven manipulation of strongly nonlinear systems.
In particular, we employ model predictive control to track a reference value in this intrinsic coordinate system. 

The proposed approach is fully data-driven, does not impose or assume any structure, except for a suitable choice of basis in which the conserved quantities are representable.
Further, we demonstrate how the effect of actuation on the eigenfunctions can be identified.
In general, identifying a parsimonious model structure through a sparsity constraint is beneficial as it prevents overfitting and has improved robustness to noise. 
However, finding a sparse vector in the null space of a matrix is a challenge, and future work is required to determine a robust algorithmic solution for identification and control using limited and noisy measurements.

To characterize and control general nonlinear systems, the KRONIC framework must be extended to approximate eigenfunctions with non-zero eigenvalue. One possible direction involves the alternating optimization of the eigenvalue and eigenfunction.
In addition, there exist several interesting generalizations of Noether's theorem for non-conservative systems.
The framework may be generalized to non-conservative systems, e.g.\ based on a variational formulation of the virtual work functional instead of action functional, which gives rise to balance laws instead of conservation laws.
Connecting these extension with the Koopman operator theory may provide a data-driven framework for the discovery and control of these balance laws. 
Finally, it will be interesting to apply this method to data from other dynamical systems with unknown symmetries and conservation laws.

\section{ACKNOWLEDGMENTS}

EK gratefully acknowledges support by the Washington Research Foundation, the Gordon and Betty Moore Foundation (Award \#2013-10-29), the Alfred P. Sloan Foundation (Award \#3835), and the University of Washington eScience Institute. 
SLB acknowledges support from the Army Research Office (ARO W911NF-17-1-0306).  
SLB and JNK acknowledge support from the Defense
Advanced Research Projects Agency (DARPA HR011-16-C-0016, PA-18-01-FP-125).


\bibliographystyle{plain}        
\bibliography{references} 

\begin{thebibliography}{10}

\bibitem{Abraham1978book}
R.~Abraham and J.~E. Marsden.
\newblock {\em Foundations of mechanics}, volume~36.
\newblock Benjamin/Cummings Publishing Company Reading, Massachusetts, 1978.

\bibitem{arnol2013mathematical}
V.~I. Arnol'd.
\newblock {\em Mathematical methods of classical mechanics}, volume~60.
\newblock Springer Science \& Business Media, 2013.

\bibitem{Bloch1996arma}
A.~M. Bloch, P.~S. Krishnaprasad, J.~E. Marsden, and R.~M. Murray.
\newblock Nonholonomic mechanical systems with symmetry.
\newblock {\em Archive for Rational Mechanics and Analysis}, 136(1):21--99,
  1996.

\bibitem{Brunton2016plosone}
S.~L. Brunton, B.~W. Brunton, J.~L. Proctor, and J.~N Kutz.
\newblock Koopman invariant subspaces and finite linear representations of
  nonlinear dynamical systems for control.
\newblock {\em PLoS ONE}, 11(2):e0150171, 2016.

\bibitem{Brunton2016pnas}
S.~L. Brunton, J.~L. Proctor, and J.~N. Kutz.
\newblock Discovering governing equations from data by sparse identification of
  nonlinear dynamical systems.
\newblock {\em Proceedings of the National Academy of Sciences},
  113(15):3932--3937, 2016.

\bibitem{camacho2013model}
E.~F. Camacho and C.~B. Alba.
\newblock {\em Model predictive control}.
\newblock Springer Science \& Business Media, 2013.

\bibitem{Chartrand2011isrnam}
R.~Chartrand.
\newblock Numerical differentiation of noisy, nonsmooth data.
\newblock {\em ISRN Appl. Math.}, 2011, 2011.

\bibitem{cimen2008ifac}
T.~Cimen.
\newblock State-dependent {R}iccati equation ({SDRE}) control: A survey.
\newblock {\em IFAC Proceedings Volumes}, 41(2):3761--3775, 2008.

\bibitem{garcia1989model}
C.~E. Garcia, D.~M. Prett, and M.~Morari.
\newblock Model predictive control: theory and practice---a survey.
\newblock {\em Automatica}, 25(3):335--348, 1989.

\bibitem{guckenheimer_holmes}
P.~Holmes and J.~Guckenheimer.
\newblock {\em Nonlinear oscillations, dynamical systems, and bifurcations of
  vector fields}, volume~42 of {\em Applied Mathematical Sciences}.
\newblock Springer-Verlag, Berlin, Heidelberg, 1983.

\bibitem{Kaiser2017arxiv}
E.~Kaiser, J.~N. Kutz, and S.~L. Brunton.
\newblock Data-driven discovery of {K}oopman eigenfunctions for control.
\newblock {\em arXiv preprint arXiv:1707.01146}, 2017.

\bibitem{koon2000heteroclinic}
W.~S. Koon, M.~W. Lo, J.~E. Marsden, and S.~D. Ross.
\newblock Heteroclinic connections between periodic orbits and resonance
  transitions in celestial mechanics.
\newblock {\em Chaos: An Interdisciplinary Journal of Nonlinear Science},
  10(2):427--469, 2000.

\bibitem{koon2008dynamical}
W.~S. Koon, M.~W. Lo, J.~E. Marsden, and S.~D. Ross.
\newblock Dynamical systems, the three-body problem and space mission design.
\newblock {\em available online. Accessed}, 21:9, 2008.

\bibitem{Koopman1931pnas}
B.~O. Koopman.
\newblock Hamiltonian systems and transformation in {H}ilbert space.
\newblock {\em Proceedings of the National Academy of Sciences},
  17(5):315--318, 1931.

\bibitem{Korda2016arxiv}
M.~Korda and I.~Mezi{\'c}.
\newblock Linear predictors for nonlinear dynamical systems: Koopman operator
  meets model predictive control.
\newblock {\em arXiv preprint arXiv:1611.03537}, 2016.

\bibitem{Kutz2016book}
J.~N. Kutz, S.~L. Brunton, B.~W. Brunton, and J.~L. Proctor.
\newblock {\em Dynamic Mode Decomposition: Data-Driven Modeling of Complex
  Systems}.
\newblock SIAM, 2016.

\bibitem{Langley1981bacon}
P.~Langley, G.~L. Bradshaw, and H.~A. Simon.
\newblock Bacon. 5: The discovery of conservation laws.
\newblock In {\em IJCAI}, volume~81, pages 121--126, 1981.

\bibitem{Marsden:MS}
J.~E. Marsden and T.~S. Ratiu.
\newblock {\em Introduction to mechanics and symmetry}.
\newblock Springer-Verlag, 2nd edition, 1999.

\bibitem{Marsden2001dmvi}
J.~E. Marsden and M.~West.
\newblock Discrete mechanics and variational integrators.
\newblock {\em Acta Numerica}, 10:357--514, 2001.

\bibitem{mezic2003cdc}
I.~Mezic.
\newblock Controllability of hamiltonian systems with drift: Action-angle
  variables and ergodic partition.
\newblock In {\em Decision and Control, 2003. Proceedings. 42nd IEEE Conference
  on}, volume~3, pages 2585--2592. IEEE, 2003.

\bibitem{Mezic2005nd}
I.~Mezi{\'c}.
\newblock Spectral properties of dynamical systems, model reduction and
  decompositions.
\newblock {\em Nonlinear Dynamics}, 41(1-3):309--325, 2005.

\bibitem{Mezic2017arxiv}
I.~Mezic.
\newblock Koopman operator spectrum and data analysis.
\newblock {\em arXiv preprint arXiv:1702.07597}, 2017.

\bibitem{Mezic2017book}
I.~Mezi\'c.
\newblock {\em Spectral operator methods in dynamical systems: Theory and
  applications}.
\newblock Draft manuscript UCSB, 2017.

\bibitem{Mezic2004physicad}
I.~Mezi{\'c} and A.~Banaszuk.
\newblock Comparison of systems with complex behavior.
\newblock {\em Physica D}, 197(1):101--133, 2004.

\bibitem{mezic1999chaos}
I.~Mezi{\'c} and S.~Wiggins.
\newblock A method for visualization of invariant sets of dynamical systems
  based on the ergodic partition.
\newblock {\em Chaos: An Interdisciplinary Journal of Nonlinear Science},
  9(1):213--218, 1999.

\bibitem{noether1918invariante}
E~Noether.
\newblock Invariante variationsprobleme nachr. d. k{\"o}nig. gesellsch. d.
  wiss. zu g{\"o}ttingen, math-phys. klasse 1918: 235-257.
\newblock {\em English Reprint: physics/0503066, http://dx. doi.
  org/10.1080/00411457108231446}, page~57, 1918.

\bibitem{Peitz2017arxiv}
S.~Peitz and S.~Klus.
\newblock Koopman operator-based model reduction for switched-system control of
  pdes.
\newblock {\em arXiv preprint arXiv:1710.06759}, 2017.

\bibitem{Qu2014}
Q.~Qu, J.~Sun, and J.~Wright.
\newblock {Finding a sparse vector in a subspace: Linear sparsity using
  alternating directions}.
\newblock In {\em Advances in Neural Information Processing Systems 27}, pages
  3401----3409, 2014.

\bibitem{Rowley2009jfm}
C.~W. Rowley, I.\ Mezi\'c, S.\ Bagheri, P.\ Schlatter, and D.S. Henningson.
\newblock Spectral analysis of nonlinear flows.
\newblock {\em J.\ Fluid Mech.}, 645:115--127, 2009.

\bibitem{Rudy2017sciadv}
S.~H. Rudy, S.~L. Brunton, J.~L. Proctor, and J.~N. Kutz.
\newblock Data-driven discovery of partial differential equations.
\newblock {\em Science Advances}, 3(e1602614), 2017.

\bibitem{Schmid2010jfm}
P.~J. Schmid.
\newblock Dynamic mode decomposition of numerical and experimental data.
\newblock {\em Journal of Fluid Mechanics}, 656:5--28, August 2010.

\bibitem{Schmidt2009science}
M.~Schmidt and H.~Lipson.
\newblock Distilling free-form natural laws from experimental data.
\newblock {\em Science}, 324(5923):81--85, 2009.

\bibitem{schulte2010discovery}
O.~Schulte and M.~S. Drew.
\newblock Discovery of conservation laws via matrix search.
\newblock In {\em International Conference on Discovery Science}, pages
  236--250. Springer, 2010.

\bibitem{Thiffeault2000physicad}
J.-L. Thiffeault and P.~J. Morrison.
\newblock Classification and casimir invariants of lie--poisson brackets.
\newblock {\em Physica D}, 136:205--244, 2000.

\bibitem{Williams2015jnls}
M.~O. Williams, I.~G. Kevrekidis, and C.~W. Rowley.
\newblock A data-driven approximation of the {K}oopman operator: extending
  dynamic mode decomposition.
\newblock {\em Journal of Nonlinear Science}, 6:1307--1346, 2015.

\bibitem{Yoshida1990pla}
H.~Yoshida.
\newblock Construction of higher order symplectic integrators.
\newblock {\em Physics letters A}, 150(5-7):262--268, 1990.

\end{thebibliography}

%
%
%
%

\end{document}